\documentclass[12pt]{aptpub}

\renewcommand{\crnotice}{\hfill{{\small\em

\/}\ {\small}}}

\usepackage{amsmath,amstext,url}
\usepackage{color}

\renewcommand{\authornames}{\underline{\footnotesize }
Junping Li}
\renewcommand{\shorttitle}
{ \underline{{Birth-Death processes with two-type catastrophes}}}

\numberwithin{equation}{section} 

\makeatletter \@addtoreset{equation}{section}

\makeatletter \@addtoreset{lemma}{section}

\makeatletter \@addtoreset{theorem}{section}

\makeatletter \@addtoreset{proposition}{section}

\makeatletter \@addtoreset{corollary}{section}

\makeatletter \@addtoreset{remark}{section}

\makeatletter \@addtoreset{definition}{section}

\makeatletter \@addtoreset{example}{section}

\begin{document}

\thispagestyle{firstpg}

\vspace*{1.5pc} \noindent \normalsize\textbf{\Large {Birth-Death processes with two-type catastrophes}} \hfill

\vspace{12pt} 
\hspace*{0.75pc}{\small\textrm{\uppercase{Junping Li}}}
\hspace{-2pt}$^{**}$, {\small\textit{Guangdong University of Science $\&$ Technology; Central
South University}}


\par
\footnote{\hspace*{-0.75pc}$^{*}\,$Postal
address: Guangdong University of Science $\&$ Technology, Dongguan, 523083, China; Central
South University, Changsha, 410083, China. E-mail:
jpli@mail.csu.edu.cn}

\par
\renewenvironment{abstract}{%
\vspace{8pt} \vspace{0.1pc} \hspace*{0.25pc}
\begin{minipage}{14cm}
\footnotesize
{\bf Abstract}\\[1ex]
\hspace*{0.5pc}} {\end{minipage}}
\begin{abstract}
  This paper concentrates on the general birth-death processes with two different types of catastrophes. The Laplace transform of transition probability function for birth-death processes with two-type catastrophes are is successfully expressed with the Laplace transform of transition probability function of the birth-death processes without catastrophe. The first effective catastrophe occurrence time is considered. The Laplace transform of its probability density function, expectation and variance are obtained.
\end{abstract}

\vspace*{12pt}
\parbox[b]{26.75pc}{{
}}
\par
\normalsize

\renewcommand{\amsprimary}[1]{
\vspace*{8pt}
\hspace*{2.25pc}
\parbox[b]{20.75pc}{\scriptsize
AMS 2000 Subject Classification: Primary 60J27 Secondary 60J35
     {\uppercase{#1}}}\par\normalsize}
\renewcommand{\ams}[2]{
\vspace*{8pt}
\parbox[b]{24.75pc}{\scriptsize
     AMS 2000 SUBJECT CLASSIFICATION: PRIMARY
     {\uppercase{#1}}\\ \phantom{
     AMS 2000
     SUBJECT CLASSIFICATION:
     }
SECONDARY
 {\uppercase{#2}}}\par\normalsize}

\ams{60J27}{60J35}

\par
\vspace{5mm}
 \setcounter{section}{1}
 \setcounter{equation}{0}
 \setcounter{theorem}{0}
 \setcounter{lemma}{0}
 \setcounter{corollary}{0}
\noindent {\large \bf 1. Introduction}
\vspace{3mm}
\par
Markov process is a very important branch of stochastic processes and has a very wide range of applications. Many research works can be referenced, such as, Anderson~\cite{A1991}, Asmussen~\cite{A2003}, Chen~\cite{CMF04} and others.
\par
The birth-death process is a very important class of Markov processes, which has been widely applied in finance, communications, population science and queueing theory. In the past few decades, there are many works on generalizing the ordinary birth-death process and make the theory of birth-death processes more and more fruitful. Recently, the stochastic models with catastrophe have aroused much research interest. For example, Chen Zhang and Liu~\cite{CZL2004}, Economou and Fakinos~\cite{EF2003}, Pakes~\cite{P1997} considered the instantaneous distribution of continuous-time Markov chains with catastrophes. Chen and Renshaw~\cite{CR1997,CR2004} analyzed the effect of catastrophes on the $M/M/1$ queuing model. Zhang and Li~\cite{ZL2015} extended these results to the $M/M/c$ queuing model with catastrophes. Li and Zhang~\cite{LZ2017} further considered the effect of catastrophes on the $M^X/M/c$ queuing model.
Di Crescenzo et al~\cite{DGNR2008} discussed the probability distribution and the relevant numerical characteristics of the first occurrence time of an effective disaster for general birth-death process with catastrophes.
Other related works can be seen from Artalejo~\cite{A2000}, Bayer and Boxma~\cite{BB1996}, Chen, Pollett, Li and Zhang~\cite{CPLZ2010}, Dudin and Karolik~\cite{DK2001}, Gelenbe~\cite{Ge1991}, Gelenbe, Glynn and Sigman~\cite{Geetal1991}, Jain and Sigman~\cite{JS1996},
\par
In this paper, we mainly consider the property of the first occurrence time of effective catastrophe for the general birth-death processes with  two-type catastrophes.
\par
We start our discussion by presenting the infinitesimal generator, i.e., the so called $q$-matrix.
\par
\begin{definition}\label{def1.1}\ Let $\{N_t:t\geq 0\}$ be a continuous-time Markov chain on state space $\mathbf{Z}_+=\{0,1,2,\cdots\}$, if its $q$-matrix $Q=(q_{ij}:i,j\in \mathbf{Z}_+)$ is by
\begin{eqnarray}\label{eq1-1}
Q=\hat{Q}+Q_d,
\end{eqnarray}
where $\hat{Q}=(\hat{q}_{ij}:i,j\in \mathbf{Z}_+)$ and $Q_d=(q^{(d)}_{ij}:i,j\in \mathbf{Z}_+)$ are given by
\begin{eqnarray}\label{eq1-2}
\hat{q}_{ij}=
\begin{cases}
\lambda_i, & \mbox{$i\geq 0,\ j=i+1$},\\
\mu_i, & \mbox{$i\geq 1,\ j=i-1$},\\
-\lambda_0,\ &\mbox{$i=j=0$},\\
-\omega_i, & \mbox{$i=j\geq 1$},\\
0,\ & \text{otherwise}.
\end{cases}
\end{eqnarray}
and
\begin{eqnarray}\label{eq1-3}
q^{(d)}_{ij}=
\begin{cases}
\beta, & \mbox{$i=0$\ or\ $i\geq 2$,\ $j=1$},\\
\alpha, & \mbox{$i\geq 1,\ j=0$},\\
-\beta,\ &\mbox{$i=j=0$},\\
-\alpha, & \mbox{$i=j=1$},\\
-\gamma, & \mbox{$i=j\geq 2$},\\
0,\ & \text{otherwise}.
\end{cases}
\end{eqnarray}
with $\alpha, \beta\geq 0, \lambda_i>0\ (i\geq 0),\ \mu_i>0\ (i\geq 1) $ and $\omega_i=\lambda_i+\mu_i\ (i\geq 1),\ \gamma=\alpha+\beta$, respectively.
\par
Then $\{N_t:t\geq 0\}$ is called a birth-death processes with two-type catastrophes. Its probability transition function is denoted by $P(t)=(p_{ij}(t):i,j\in \mathbf{Z}_+)$ and the corresponding resolvent is denoted by $\Pi(\lambda)=(\pi_{j,n}(\lambda):j,n\in \mathbf{Z}_+)$.
\end{definition}
\par
\begin{remark}\label{re1.1}\ By Definition~\ref{def1.1}, $\alpha$ and $\beta$ describe the rates of catastrophes. We called them $\alpha$-catastrophe and $\beta$-catastrophe, respectively. That is, $\alpha$-catastrophe kills all the individuals in the system, while $\beta$-catastrophe partially kills the individuals in the system with only one individual left. If $\alpha=\beta=0$, i.e., there is no catastrophe, then $\{N_t:t\geq 0\}$ degenerates into an ordinary birth-death process, which is denoted by $\{\hat{N}(t):t\geq 0\}$, its $q$-matrix is denoted by $\hat{Q}$.
The probability transition function of $\{\hat{N}_t:t\geq 0\}$ is denoted by $\hat{P}(t)=(\hat{p}_{ij}(t):i,j\in \mathbf{Z}_+)$ and the corresponding resolvent is denoted by $\hat{\Pi}(\lambda)=(\hat{\pi}_{j,n}(\lambda):j,n\in \mathbf{Z}_+)$.
\end{remark}

\par
\vspace{5mm}
 \setcounter{section}{2}
 \setcounter{equation}{0}
 \setcounter{theorem}{0}
 \setcounter{lemma}{0}
 \setcounter{definition}{0}
 \setcounter{corollary}{0}
\noindent {\large \bf 2. Probability transition function}
 \vspace{3mm}
\par
From Definition~\ref{def1.1}, we see that a catastrophe may reduce the system state to $0$ or $1$. However, since natural death rate $\mu_1,\ \mu_2>0$, when the system state transfer to $0$ from $1$ or transfer to $1$ from $2$, it is difficult to distinguish whether it was a catastrophe or a natural death. Therefore, it is important to discuss such effective catastrophe. For this purpose, we first construct the relationship of $P(t)$ and $\hat{P}(t)$ (or equivalently, $\Pi(\lambda)$ and $\hat{\Pi}(\lambda)$).
\par
\begin{lemma}\label{le2.1}\ {\rm(i)}\ $P(t)=(p_{j,n}(t):j,n\in\mathbf{Z}_+)$ satisfies the following Kolmogorov forward equations: for any $j,n\in \mathbf{Z}_+$ and $t\geq 0$,
\begin{eqnarray}\label{eq2-1}
\begin{cases}
p'_{j,0}(t)=-(\lambda_0+\gamma)p_{j,0}(t)+\mu_1 p_{j,1}(t)+\alpha,\\
p'_{j,1}(t)=\lambda_0p_{j,0}(t)-(\omega_1+\gamma)p_{j,1}(t)
+\mu_2p_{j,2}(t)+\beta,\\
p'_{j,n}(t)=\lambda_{n-1}p_{j,n-1}(t)-(\omega_n+\gamma)p_{j,n}(t)
+\mu_{n+1}p_{j,n+1}(t),\ \ n\geq 2,
\end{cases}
\end{eqnarray}
or equivalently, in the resolvent version,
\begin{eqnarray}\label{eq2-2}
\begin{cases}
(\lambda+\lambda_0+\gamma)\pi_{j,0}(\lambda)-\delta_{j,0}=\mu_1 \pi_{j,1}(\lambda)+\frac{\alpha}{\lambda},\\
(\lambda+\omega_1+\gamma) \pi_{j,1}(\lambda)-\delta_{j,1}=\lambda_0\pi_{j,0}(\lambda)+\mu_2\pi_{j,2}(\lambda)+\frac{\beta}{\lambda},\\
(\lambda+\omega_n+\gamma)\pi_{j,n}(\lambda)-\delta_{j,n}=\lambda_{n-1}\pi_{j,n-1}(\lambda)
+\mu_{n+1}\pi_{j,n+1}(\lambda),\ \ n\geq 2.
\end{cases}
\end{eqnarray}
\par
{\rm(ii)}\ $\hat{P}(t)=(\hat{p}_{j,n}(t):j,n\in\mathbf{Z}_+)$ satisfies the following Kolmogorov forward equations: for any $j,n\in \mathbf{Z}_+$ and $t\geq 0$,
\begin{eqnarray*}
\begin{cases}
\hat{p}'_{j,0}(t)=-\lambda_0\hat{p}_{j,0}(t)+\mu_1 \hat{p}_{j,1}(t),\\
\hat{p}'_{j,n}(t)=\lambda_{n-1}\hat{p}_{j,n-1}(t)-(\lambda_n+\mu_n)\hat{p}_{j,n}(t)
+\mu_{n+1}\hat{p}_{j,n+1}(t),\ \ n\geq 1.
\end{cases}
\end{eqnarray*}
or equivalently, in the resolvent version,
\begin{eqnarray*}
\begin{cases}
(\lambda+\lambda_0)\hat{\pi}_{j,0}(\lambda)-\delta_{j,0}=\mu_1 \hat{\pi}_{j,1}(\lambda),\\
(\lambda+\lambda_n+\mu_n)\hat{\pi}_{j,n}(\lambda)-\delta_{j,n}=\lambda_{n-1}\hat{\pi}_{j,n-1}(\lambda)
+\mu_{n+1}\hat{\pi}_{j,n+1}(\lambda),\ \ n\geq 1.
\end{cases}
\end{eqnarray*}
\end{lemma}
\par
\begin{proof}
(i)\ By Kolmogorov forward equations and the honesty of $P(t)$, we know that
\begin{eqnarray*}
p'_{j,0}(t)&=&-(\lambda_0+\beta)p_{j,0}(t)+(\mu_1+\alpha)p_{j,1}(t)+
\sum\limits_{k=2}^{\infty}\alpha p_{j,k}(t)\\
&=&-(\lambda_0+\beta)p_{j,0}(t)+\mu_1p_{j,1}(t)+
\sum\limits_{k=1}^{\infty}\alpha p_{j,k}(t)\\
&=&-(\lambda_0+\beta)p_{j,0}(t)+\mu_1p_{j,1}(t)+\alpha(1-p_{j,0}(t))\\
&=&-(\lambda_0+\gamma)p_{j,0}(t)+\mu_1p_{j,1}(t)+\alpha.
\end{eqnarray*}
and
\begin{eqnarray*}
p'_{j,1}(t)&=&(\lambda_0+\beta)p_{j,0}(t)-(\lambda_1+\mu_1+\alpha)p_{j,1}(t)+
(\mu_2+\beta)p_{j,2}(t)+\sum\limits_{k=3}^{\infty}\beta p_{j,k}(t)\\
&=&\lambda_0p_{j,0}(t)-(\omega_1+\alpha)p_{j,1}(t)+
\mu_2p_{j,2}(t)+\beta(1-p_{j,1}(t))\\
&=&\lambda_0p_{j,0}(t)-(\omega_1+\gamma)p_{j,1}(t)
+\mu_2p_{j,2}(t)+\beta.
\end{eqnarray*}
The other equalities of (i) and (ii) follow directly from Kolmogorov forward equations and Laplace transform. The proof is complete. \hfill $\Box$
\end{proof}
\par
The following theorem plays an important role in the later discussion, it reveals the relationship of $P(t)$ and $\hat{P}(t)$ (or equivalently, $\Pi(\lambda)$ and $\hat{\Pi}(\lambda)$).
\par
\begin{theorem}\label{th2.1}
For any $j,n\in\mathbf{Z}_+$, we have
\par
\begin{eqnarray}\label{eq2-3}
p_{j,n}(t)&=&e^{-\gamma t}\hat{p}_{j,n}(t)+\alpha\int_0^t e^{-\gamma s}\hat{p}_{0,n}(s)ds+\beta\int_0^te^{-\gamma s}\hat{p}_{1,n}(s)ds
\end{eqnarray}
or equivalently in resolvent version,
\begin{eqnarray}\label{eq2-4}
\pi_{j,n}(\lambda)&=&\hat{\pi}_{j,n}(\lambda+\gamma)+
\frac{1}{\lambda}\cdot [
\alpha\hat{\pi}_{0,n}(\lambda+\gamma)+
\beta\hat{\pi}_{1,n}(\lambda+\gamma)]
\end{eqnarray}
\end{theorem}
\par
\begin{proof}\ We first assume $\alpha=0$. The corresponding process is denoted by $\tilde{N}_t$ and its probability transition function is denoted by $\tilde{P}(t)=(\tilde{p}_{j,n}(t):j,n\in \mathbf{Z}_+)$. Denote $\{A_t:t\geq 0\}=\{\hat{N}_t:t\geq 0\}$. Let $\{K_t:t\geq 0\}$ be a Poisson process with parameter $\beta$, which is independent of $\{A_t:t\geq 0\}$, note that $\{K_t:t\geq 0\}$ can be viewed as a catastrophe flow. Let $l(t)$ be the time until the first catastrophe before time $t$. Then $l(t)$ has the truncated exponential law
\begin{eqnarray*}
P(l(t)\leq u)=1-e^{-\beta u}I_{[0,t)}(u).
\end{eqnarray*}
\par
 Denote $\{A^{(0)}_t:t\geq 0\}:=\{A_t:t\geq 0\}$. Let $\{A^{(n)}_t:t\geq 0\}_{n\geq 1}$ be an independent sequence copies of $\{A^{(0)}_t:t\geq 0\}$ but with $A^{(n)}_0=1$. Define $\{R_t:t\geq 0\}$ by
\begin{eqnarray*}
R_t=A^{(K_t)}_{l(t)},\quad t\geq 0.
\end{eqnarray*}
Then, $\{R_t:t\geq 0\}$ is a continuous-time Markov chain, it evolves like  $A^{(0)}_t$, at the first catastrophe time, it jumps to state $1$, and then evolves like $A^{(1)}_t$, at the next catastrophe time, it jumps to state $1$ again, and so on. Let $\bar{P}(t)=(\bar{p}_{jn}(t):j,n\in \mathbf{Z}_+))$ be the probability transition function of $\{R_t:t\geq 0\}$. Then
\begin{eqnarray*}
\bar{p}_{jn}(t)=P(R_t=n|R_0=j)=P_j(R_t=n)=E_j[I_{\{n\}}(R_t)]
=E_j[E_j[I_{\{n\}}(A^{(K_t)}_{l(t)})|K_t,l(t)]],
\end{eqnarray*}
where $P_j=P(\cdot |R_0=j)$ and $E_j$ is the mathematical expectation under $P_j$. Denote $G(K_t,l(t)):=E_j[I_{\{n\}}(A^{(K_t)}_{l(t)})|K_t,l(t)]$ for a moment. Then the above equality equals to
\begin{eqnarray*}
&&E_j[G(K_t,l(t))]\\
&=& E_j[E_j[G(K_t,l(t))|l(t)]]\\
&=&P_j(l(t)=t)E_j[G(K_t,l(t))|l(t)=t]+\beta\xi\int_0^t e^{-\beta  s}E_j[G(K_t,l(t))|l(t)=s]ds.
\end{eqnarray*}
Since $l(t)=t\Leftrightarrow K_t=0$ and $R_0=j\Leftrightarrow A_0=j$, we have
\begin{eqnarray*}
P_j(l(t)=t)=P_j(K_t=0)=e^{-\beta t}
\end{eqnarray*}
and
\begin{eqnarray*}
E_j[G(K_t,l(t))|l(t)=t]=E_j[I_{\{n\}}(A^{(0)}_{t})]
=E_j[I_{\{n\}}(A_{t})]=\hat{p}_{jn}(t).
\end{eqnarray*}
If $s<t$, then
\begin{eqnarray*}
&&E_j[G(K_t,l(t))|l(t)=s]\\
&=&\sum\limits_{k=1}^{\infty}P_j(K_t=k|l(t)=s)G(k,s)\\
&=&\sum\limits_{k=1}^{\infty}P_j(K_t=k|l(t)=s)
E_j[I_{\{n\}}(A^{(K_t)}_{l(t)})|K_t=k,l(t)=s]\\
&=&\sum\limits_{k=1}^{\infty}P_j(K_t=k|l(t)=s)
E_j[I_{\{n\}}(A^{(k)}_{s})]\\
&=&\sum\limits_{k=1}^{\infty}P_j(K_t=k|l(t)=s)
E[I_{\{n\}}(A^{(k)}_{s})|A_0=j,A^{(k)}_0=1]\\
&=&\sum\limits_{k=1}^{\infty}P_j(K_t=k|l(t)=s)
P(A^{(k)}_{s}=n|A_0=j,A^{(k)}_0=1]\\
&=&\sum\limits_{k=1}^{\infty}P_j(K_t=k|l(t)=s)
P(A^{(k)}_{s}=n|A^{(k)}_0=1]\\
&=&\sum\limits_{k=1}^{\infty}P_j(K_t=k|l(t)=s)
P(A_{s}=n|A_0=1]\\
&=&\hat{p}_{1,n}(s).\\
\end{eqnarray*}
Therefore,
\begin{eqnarray*}
\bar{p}_{j,n}(t)=e^{-\beta t}\hat{p}_{j,n}(t)+\beta\xi\int_0^t e^{-\beta s}\hat{p}_{1,n}(s)ds.
\end{eqnarray*}
It is easy to check that $\bar{p}'_{j,n}(0)=\tilde{p}'_{j,n}(0)$. This implies that $R_t$ and $\tilde{N}_t$ are same in sense of distribution. Hence,
\begin{eqnarray}\label{eq2-5}
\tilde{p}_{j,n}(t)=e^{-\beta t}\hat{p}_{j,n}(t)+\beta\int_0^t e^{-\beta s}\hat{p}_{1,n}(s)ds.
\end{eqnarray}
\par
Now consider the general case $\alpha>0$. Denote $\{\tilde{A}_t:t\geq 0\}:=\{\tilde{N}_t:t\geq 0\}$. Let $\{\tilde{K}_t:t\geq 0\}$ be a Poisson process with parameter $\alpha\xi$, which is independent of $\{\tilde{A}_t:t\geq 0\}$. $\{\tilde{K}_t:t\geq 0\}$ can be viewed as a catastrophe flow with parameter $\alpha$. Let $\tilde{l}(t)$ be the time until the first catastrophe before
time $t$. Then $l(t)$ has the truncated exponential law
\begin{eqnarray*}
P(\tilde{l}(t)\leq u)=1-e^{-\alpha u}I_{[0,t)}(u).
\end{eqnarray*}
\par
Denote $\{\tilde{A}^{(0)}_t:t\geq 0\}:=\{\tilde{A}_t:t\geq 0\}$. Let $\{\tilde{A}^{(n)}_t:t\geq 0\}_{n\geq 1}$ be an independent sequence copies of $\{\tilde{A}^{(0)}_t:t\geq 0\}$ but with $\tilde{A}^{(n)}_0=0\ (n\geq 1)$. Define $\{\tilde{R}_t:t\geq 0\}$ by
\begin{eqnarray*}
\tilde{R}_t=\tilde{A}^{(\tilde{K}_t)}_{\tilde{l}(t)},\quad t\geq 0.
\end{eqnarray*}
Let $\check{P}(t)=(\check{p}_{j,n}(t):j,n\in\mathbf{Z}_+)$ be the probability transition function of $\{\tilde{R}_t:t\geq 0\}$. By a similar argument as above, we know that
\begin{eqnarray*}
\check{p}_{j,n}(t)=e^{-\alpha t}\bar{p}_{j,n}(t)+\alpha\int_0^t e^{-\alpha s}\bar{p}_{0,n}(s)ds.
\end{eqnarray*}
By (\ref{eq2-5})
\begin{eqnarray*}
\check{p}_{j,n}(t)&=&e^{-\alpha t}[e^{-\beta t}\hat{p}_{j,n}(t)+\beta\int_0^t e^{-\beta s}\hat{p}_{1,n}(s)ds]\\
&&+\alpha\int_0^t e^{-\alpha s}[e^{-\beta s}\hat{p}_{0,n}(s)+\beta\int_0^s e^{-\beta u}\hat{p}_{1,n}(u)du]ds\\
&=&e^{-(\alpha+\beta)t}\hat{p}_{j,n}(t)+\alpha\int_0^t e^{-(\alpha+\beta) s}\hat{p}_{0,n}(s)ds+\beta\int_0^te^{-(\alpha+\beta) s}\hat{p}_{1,n}(s)ds.\\
\end{eqnarray*}
It is easy to check that $\check{p}'_{j,n}(0)=p'_{j,n}(0)$. This implies that $\tilde{R}_t$ and $N_t$ are same in sense of distribution. Hence,
\begin{eqnarray*}
p_{j,n}(t)=e^{-(\alpha+\beta)t}\hat{p}_{j,n}(t)+\alpha\int_0^t e^{-(\alpha+\beta) s}\hat{p}_{0,n}(s)ds+\beta\int_0^te^{-(\alpha+\beta) s}\hat{p}_{1,n}(s)ds.\\
\end{eqnarray*}
(\ref{eq2-3}) is proved. Taking Laplace transform on (\ref{eq2-3}) implies (\ref{eq2-4}). The proof is complete.
\hfill $\Box$
\end{proof}

\par
\vspace{5mm}
 \setcounter{section}{3}
 \setcounter{equation}{0}
 \setcounter{theorem}{0}
 \setcounter{lemma}{0}
 \setcounter{definition}{0}
 \setcounter{corollary}{0}
\noindent {\large \bf 3.  The first occurrence time of effective catastrophe}
\vspace{3mm}
\par
We now consider the first effective catastrophe of $\{N_t:t\geq 0\}$. Let $C_j$ is the first occurrence time of effective catastrophe for $\{N_t:t\geq 0\}$ starting from state $j$. The probability density function of $C_j$ is denoted by $d_j(t)$.
Let $C_{j,0}$ and $C_{j,1}$ be the first occurrence time of effective $\alpha$-catastrophe and effective $\beta$-catastrophe, respectively. It is obvious that $C_j=C_{j,0}\wedge C_{j,1}$.
\par
The property of $C_{j,0}$ or $C_{j,1}$ can be similarly discussed as in Di Crescenzo et al~\cite{DGNR2008}. In this paper, we mainly consider the property of $C_j$ and the probabilities $P(C_j\leq t, C_{j,0}<C_{j,1})$ and $P(C_j\leq t, C_{j,1}<C_{j,0})$. For this purpose, we construct a new process $\{M_t:t\geq 0\}$ such that $\{M_t:t\geq 0\}$ coincides with $\{N_t:t\geq 0\}$ until the occurrence of catastrophe, but $\{M_t:t\geq 0\}$ enter into an absorbing state $-1$ if the first effective catastrophe is $\beta$-type and enter into another absorbing state $-2$ if the first effective catastrophe is $\alpha$-type. Therefore the state space of $\{M_t:t\geq 0\}$ is $\mathbf{S}=\{-2,-1,0,1,\cdots\}$ and its $q$-matrix $\tilde{Q}=(\tilde{q}_{jn}:j,n\in \mathbf{S})$ is given by
\begin{eqnarray*}
\tilde{q}_{ij}=
\begin{cases}
\lambda_i, & \mbox{$i\geq 0,\ j=i+1$},\\
\mu_i, & \mbox{$i\geq1,\ j=i-1$},\\
\alpha, & \mbox{$i\geq 1,\ j=-2$},\\
\beta, & \mbox{$i=0,\ j=-1$},\\
\beta, & \mbox{$i\geq 2,\ j=-1$},\\
-(\lambda_0+\beta),\ &\mbox{$i=j=0$},\\
-(\omega_1+\alpha), & \mbox{$i=j=1$},\\
-(\omega_i+\gamma), & \mbox{$i=j\geq 2$},\\
0,\ & \text{otherwise}.
\end{cases}
\end{eqnarray*}
\par
Let $H(t)=(h_{j,n}(t):j,n\in\mathbf{S})$ and $\Phi(\lambda)=(\phi_{j,n}(\lambda):j,n\in\mathbf{S})$ be the $\tilde{Q}$-transition function and $\tilde{Q}$-resolvent.
\par
\begin{lemma}\label{le3.1}
For any $j\geq 0$, we have
\begin{eqnarray}\label{eq3-1}
\begin{cases}
h'_{j,-2}(t)=\alpha (1-h_{j,-2}(t)-h_{j,-1}(t)-h_{j,0}(t)),\\
h'_{j,-1}(t)=\beta (1-h_{j,-2}(t)-h_{j,-1}(t)-h_{j,1}(t)),\\
h'_{j,0}(t)=-(\lambda_0+\beta)h_{j,0}(t)+\mu_1h_{j,1}(t),\\
h'_{j,1}(t)=\lambda_0h_{j,0}(t)-(\omega_1+\alpha)h_{j,1}(t)
+\mu_2h_{j,2}(t),\\
h'_{j,n}(t)=\lambda_{n-1}h_{j,n-1}(t)-(\omega_n+\gamma)h_{j,n}(t)
+\mu_{n+1}h_{j,n+1}(t),\quad n\geq 2,\\
\end{cases}
\end{eqnarray}
or equivalently, in resolvent version,
\begin{eqnarray}\label{eq3-2}
\begin{cases}
\lambda \phi_{j,-2}(\lambda)=\alpha (\frac{1}{\lambda}-\phi_{j,-2}(\lambda)-\phi_{j,-1}(\lambda)-\phi_{j,0}(\lambda)),\\
\lambda\phi_{j,-1}(\lambda)=\beta (\frac{1}{\lambda}-\phi_{j,-2}(\lambda)-\phi_{j,-1}(\lambda)-\phi_{j,1}(\lambda)),\\
(\lambda+\lambda_0+\beta)\phi_{j,0}(\lambda)-\delta_{j,0}
=\mu_1\phi_{j,1}(\lambda),\\
(\lambda+\omega_1+\alpha)\phi_{j,1}(\lambda)-\delta_{j,1}
=\lambda_0\phi_{j,0}(\lambda)
+\mu_2\phi_{j,2}(\lambda),\\
(\lambda+\omega_n+\gamma)\phi_{j,n}(\lambda)-\delta_{j,n}
=\lambda_{n-1}\phi_{j,n-1}(\lambda)
+\mu_{n+1}\phi_{j,n+1}(\lambda),\quad n\geq 2.\\
\end{cases}
\end{eqnarray}
\end{lemma}
\par
\begin{proof}\ By Kolmogorov forward equation,
\begin{eqnarray*}
h'_{j,-2}(t)&=&\sum\limits_{k=1}^{\infty}\alpha h_{j,k}(t)\\
&=&\alpha(1-h_{j,-2}(t)-h_{j,-1}(t)-h_{j,0}(t)).
\end{eqnarray*}
\begin{eqnarray*}
h'_{j,-1}(t)&=&\beta h_{j,0}(t)+\sum\limits_{k=2}^{\infty}\beta h_{j,k}(t)\\
&=&\beta(1-h_{j,-2}(t)-h_{j,-1}(t)-h_{j,1}(t)).
\end{eqnarray*}
The other equalities of (\ref{eq3-1}) follow directly from Kolmogorov forward equations and (\ref{eq3-2}) follows from the Laplace transform of (\ref{eq3-1}). The proof is complete. \hfill $\Box$
\end{proof}
\par
We now investigate the relationship of $\Phi(\lambda)$ and $\Pi(\lambda)$. For this purpose, define
\begin{eqnarray}\label{eq3-3}
A_{ij}(\lambda)=1-\lambda\pi_{i,j}(\lambda),\quad i,j\geq 0
\end{eqnarray}
and
\begin{eqnarray}\label{eq3-4}
H(\lambda)=\lambda^{-1}\{[\lambda+\alpha A_{00}(\lambda)][\lambda+\beta A_{11}(\lambda)]
-\alpha\beta A_{10}(\lambda)A_{01}(\lambda)\}.
\end{eqnarray}
\par
\begin{theorem}\label{th3.1}
Let $\Phi(\lambda)=(\phi_{j,n}(\lambda):j,n\in\mathbf{S})$ be the $\tilde{Q}$-resolvent. Then
\begin{eqnarray}\label{eq3-5}
\phi_{0,n}(\lambda)=\frac{(\lambda+\beta A_{11}(\lambda))\pi_{0,n}(\lambda)-\beta A_{01}(\lambda)\pi_{1,n}(\lambda)}
{H(\lambda)},\quad n\geq 0,
\end{eqnarray}
\begin{eqnarray}\label{eq3-6}
\phi_{1,n}(\lambda)=\frac{-\alpha A_{10}(\lambda)\pi_{0,n}(\lambda)+(\lambda+\alpha A_{00}(\lambda))\pi_{1,n}(\lambda)}
{H(\lambda)},\quad n\geq 0
\end{eqnarray}
and
\begin{eqnarray}\label{eq3-7}
\phi_{j,n}(\lambda)=\pi_{j,n}(\lambda)+F_j(\lambda)\pi_{0,n}(\lambda)+G_j(\lambda)\pi_{1,n}(\lambda),\quad j\geq 2,\ n\geq 0,
\end{eqnarray}
where
\begin{eqnarray}\label{eq3-8}
F_j(\lambda)=\frac{\alpha\beta A_{10}(\lambda)A_{j1}(\lambda)
-\alpha(\lambda+\beta A_{11}(\lambda))A_{j0}(\lambda)
}{\lambda H(\lambda)}
\end{eqnarray}
and
\begin{eqnarray}\label{eq3-9}
G_j(\lambda)=\frac{\alpha\beta A_{01}(\lambda)A_{j0}(\lambda)
-\beta(\lambda+\alpha A_{00}(\lambda))A_{j1}(\lambda)
}{\lambda H(\lambda)}
\end{eqnarray}
with $(\pi_{j,n}(\lambda):j,n\geq 0)$ being given by $(\ref{eq2-4})$.
\end{theorem}
\par
\begin{proof}\ By (\ref{eq3-2}) with $j=0,1$,
\begin{eqnarray}\label{eq3-10}
\begin{cases}
(\lambda+\lambda_0+\beta)\phi_{0,0}(\lambda)-1=\mu_1\phi_{0,1}(\lambda),\\
(\lambda+\omega_1+\alpha)\phi_{0,1}(\lambda)=\lambda_0\phi_{0,0}(\lambda)
+\mu_2\phi_{0,2}(\lambda),\\
(\lambda+\omega_n+\gamma)\phi_{0,n}(\lambda)=\lambda_{n-1}\phi_{0,n-1}(\lambda)
+\mu_{n+1}\phi_{0,n+1}(\lambda),\quad n\geq 2,\\
\end{cases}
\end{eqnarray}
\begin{eqnarray}\label{eq3-11}
\begin{cases}
(\lambda+\lambda_0+\beta)\phi_{1,0}(\lambda)=\mu_1\phi_{1,1}(\lambda),\\
(\lambda+\omega_1+\alpha)\phi_{1,1}(\lambda)-1=\lambda_0\phi_{1,0}(\lambda)
+\mu_2\phi_{1,2}(\lambda),\\
(\lambda+\omega_n+\gamma)\phi_{1,n}(\lambda)=\lambda_{n-1}\phi_{1,n-1}(\lambda)
+\mu_{n+1}\phi_{1,n+1}(\lambda),\quad n\geq 2\\
\end{cases}
\end{eqnarray}
and by (\ref{eq2-2}) with $j=0,1$,
\begin{eqnarray}\label{eq3-12}
\begin{cases}
(\lambda+\lambda_0+\gamma)\pi_{0,0}(\lambda)-1=\mu_1 \pi_{0,1}(\lambda)+\frac{\alpha}{\lambda},\\
(\lambda+\omega_1+\gamma) \pi_{0,1}(\lambda)=\lambda_0\pi_{0,0}(\lambda)+\mu_2\pi_{0,2}(\lambda)+\frac{\beta}{\lambda},\\
(\lambda+\omega_n+\gamma)\pi_{0,n}(\lambda)=\lambda_{n-1}\pi_{0,n-1}(\lambda)
+\mu_{n+1}\pi_{0,n+1}(\lambda),\ \ n\geq 2.
\end{cases}
\end{eqnarray}
\begin{eqnarray}\label{eq3-13}
\begin{cases}
(\lambda+\lambda_0+\gamma)\pi_{1,0}(\lambda)=\mu_1 \pi_{1,1}(\lambda)+\frac{\alpha}{\lambda},\\
(\lambda+\omega_1+\gamma) \pi_{1,1}(\lambda)-1=\lambda_0\pi_{1,0}(\lambda)+\mu_2\pi_{1,2}(\lambda)+\frac{\beta}{\lambda},\\
(\lambda+\omega_n+\gamma)\pi_{1,n}(\lambda)=\lambda_{n-1}\pi_{1,n-1}(\lambda)
+\mu_{n+1}\pi_{1,n+1}(\lambda),\ \ n\geq 2.
\end{cases}
\end{eqnarray}
\par
Let
\begin{eqnarray}\label{eq3-14}
\phi_{0,n}(\lambda)=A(\lambda)\pi_{0,n}(\lambda)+B(\lambda)\pi_{1,n}(\lambda),\quad n\geq 0.
\end{eqnarray}
Substitute (\ref{eq3-14}) into (\ref{eq3-10}) and use (\ref{eq3-12}), we have
\begin{eqnarray}\label{eq3-15}
\begin{cases}
(\lambda+\alpha A_{00}(\lambda))A(\lambda)
+\alpha A_{10}(\lambda)B(\lambda)=\lambda\\
\beta A_{01}(\lambda)A(\lambda)
+(\lambda+\beta A_{11}(\lambda))B(\lambda)=0.\\
\end{cases}
\end{eqnarray}
Indeed, by the first equality of (\ref{eq3-10}),
\begin{eqnarray*}
(\lambda+\lambda_0+\beta)[A(\lambda)\pi_{0,0}(\lambda)+B(\lambda)\pi_{1,0}
(\lambda)]-1
=\mu_1[A(\lambda)\pi_{0,1}(\lambda)+B(\lambda)\pi_{1,1}
(\lambda)]
\end{eqnarray*}
i.e.,
\begin{eqnarray*}
A(\lambda)[(\lambda+\lambda_0+\beta)\pi_{0,0}(\lambda)-
\mu_1\pi_{0,1}(\lambda)]+B(\lambda)
[(\lambda+\lambda_0+\beta)\pi_{1,0}(\lambda)-\mu_1\pi_{1,1}(\lambda)]=1
\end{eqnarray*}
It follows from the first equality of (\ref{eq3-12}) and the first equality of (\ref{eq3-13}) that
\begin{eqnarray*}
(\lambda+\alpha A_{00}(\lambda))A(\lambda)+
\alpha A_{10}(\lambda)B(\lambda)=\lambda.
\end{eqnarray*}
By the second equality of (\ref{eq3-10}),
\begin{eqnarray*}
&&(\lambda+\omega_1+\alpha)[A(\lambda)\pi_{0,1}(\lambda)+B(\lambda)\pi_{1,1}
(\lambda)]\\
&=&\lambda_0[A(\lambda)\pi_{0,0}(\lambda)+B(\lambda)\pi_{1,0}
(\lambda)]+\mu_2[A(\lambda)\pi_{0,2}(\lambda)+B(\lambda)\pi_{1,2}
(\lambda)]
\end{eqnarray*}
i.e.,
\begin{eqnarray*}
&&A(\lambda)[(\lambda+\omega_1+\alpha)\pi_{0,1}(\lambda)-\lambda_0\pi_{0,0}
(\lambda)-\mu_2\pi_{0,2}(\lambda)]\\
&&+B(\lambda)
[(\lambda+\omega_1+\alpha)\pi_{1,1}(\lambda)-\lambda_0\pi_{1,0}(\lambda)-
\mu_2\pi_{1,2}(\lambda)]=0
\end{eqnarray*}
It follows from the second equality of (\ref{eq3-12}) and the second equality of (\ref{eq3-13}) that
\begin{eqnarray*}
\beta A_{01}(\lambda)A(\lambda)+
(\lambda+\beta A_{11}(\lambda))B(\lambda)=0.
\end{eqnarray*}
Therefore, (\ref{eq3-15}) holds. It follows from (\ref{eq3-15}) that
\begin{eqnarray*}
A(\lambda)=\frac{\lambda+\beta A_{11}(\lambda)}{H(\lambda)}\quad and \quad B(\lambda)=\frac{-\beta A_{01}(\lambda)}
{H(\lambda)}.
\end{eqnarray*}
The other equalities of (\ref{eq3-10}) also hold.
\par
Let
\begin{eqnarray}\label{eq3-16}
\phi_{1,n}(\lambda)=C(\lambda)\pi_{0,n}(\lambda)+D(\lambda)\pi_{1,n}(\lambda),\quad n\geq 0.
\end{eqnarray}
Substitute (\ref{eq3-16}) into (\ref{eq3-11}) and use (\ref{eq3-13}), we have
\begin{eqnarray}\label{eq3-17}
\begin{cases}
(\lambda+\alpha-\alpha\lambda\pi_{0,0}(\lambda))C(\lambda)
+\alpha(1-\lambda\pi_{1,0}(\lambda))D(\lambda)=0\\
\beta(1-\lambda\pi_{0,1}(\lambda))C(\lambda)
+(\lambda+\beta-\beta\lambda\pi_{1,1}(\lambda))D(\lambda)=\lambda.\\
\end{cases}
\end{eqnarray}
Indeed, by the second equality of (\ref{eq3-11}),
\begin{eqnarray*}
&&(\lambda+\omega_1+\alpha)[C(\lambda)\pi_{0,1}(\lambda)+D(\lambda)\pi_{1,1}
(\lambda)]-1\\
&=&\lambda_0[C(\lambda)\pi_{0,0}(\lambda)+D(\lambda)\pi_{1,0}
(\lambda)]+\mu_2[C(\lambda)\pi_{0,2}(\lambda)+D(\lambda)\pi_{1,2}
(\lambda)]
\end{eqnarray*}
i.e.,
\begin{eqnarray*}
&&[(\lambda+\omega_1+\alpha)\pi_{0,1}(\lambda)
-\lambda_0\pi_{0,0}(\lambda)-
\mu_2\pi_{0,2}(\lambda)]C(\lambda)\\
&&+
[(\lambda+\omega_1+\alpha)\pi_{1,1}(\lambda)-\lambda_0\pi_{1,0}
(\lambda)-\mu_2\pi_{1,2}(\lambda)]D(\lambda)=1
\end{eqnarray*}
It follows from the second equality of (\ref{eq3-12}) and the second equality of (\ref{eq3-13}) that
\begin{eqnarray*}
\beta A_{01}(\lambda)C(\lambda)+
(\lambda+\beta A_{11}(\lambda))D(\lambda)=\lambda.
\end{eqnarray*}
By the first equality of (\ref{eq3-11}),
\begin{eqnarray*}
(\lambda+\lambda_0+\beta)[C(\lambda)\pi_{0,0}(\lambda)+D(\lambda)\pi_{1,0}
(\lambda)]=\mu_1[C(\lambda)\pi_{0,1}(\lambda)+D(\lambda)\pi_{1,1}
(\lambda)]
\end{eqnarray*}
i.e.,
\begin{eqnarray*}
[(\lambda+\lambda_0+\beta)\pi_{0,0}(\lambda)-\mu_1\pi_{0,1}(\lambda)]C(\lambda)+B(\lambda)
[(\lambda+\lambda_0+\beta)\pi_{1,0}(\lambda)-
\mu_1\pi_{1,1}(\lambda)]=0
\end{eqnarray*}
It follows from the first equality of (\ref{eq3-12}) and the first equality of (\ref{eq3-13}) that
\begin{eqnarray*}
(\lambda+\alpha A_{00}(\lambda))C(\lambda)+
\alpha A_{10}(\lambda))D(\lambda)=0.
\end{eqnarray*}
Therefore, (\ref{eq3-17}) holds. It follows from (\ref{eq3-17}) that
\begin{eqnarray*}
C(\lambda)=\frac{-\alpha A_{10}(\lambda)}{H(\lambda)}\quad and \quad
D(\lambda)=\frac{\lambda+\alpha A_{00}(\lambda)}
{H(\lambda)}.
\end{eqnarray*}
The other equalities of (\ref{eq3-11}) also hold.
\par
By (\ref{eq3-2}) with $j\geq 2$,
\begin{eqnarray}\label{eq3-18}
\begin{cases}
(\lambda+\lambda_0+\beta)\phi_{j,0}(\lambda)=\mu_1\phi_{j,1}(\lambda),\\
(\lambda+\omega_1+\alpha)\phi_{j,1}(\lambda)=\lambda_0\phi_{j,0}(\lambda)
+\mu_2\phi_{j,2}(\lambda),\\
(\lambda+\omega_n+\gamma)\phi_{j,n}(\lambda)-\delta_{j,n}=\lambda_{n-1}\phi_{j,n-1}(\lambda)
+\mu_{n+1}\phi_{j,n+1}(\lambda),\quad n\geq 2\\
\end{cases}
\end{eqnarray}
and by (\ref{eq2-2}) with $j\geq 2$,
\begin{eqnarray}\label{eq3-19}
\begin{cases}
(\lambda+\lambda_0+\gamma)\pi_{j,0}(\lambda)=\mu_1 \pi_{j,1}(\lambda)+\frac{\alpha}{\lambda},\\
(\lambda+\omega_1+\gamma) \pi_{j,1}(\lambda)=\lambda_0\pi_{j,0}(\lambda)+\mu_2\pi_{j,2}(\lambda)+\frac{\beta}{\lambda},\\
(\lambda+\omega_n+\gamma)\pi_{j,n}(\lambda)-\delta_{j,n}=\lambda_{n-1}\pi_{j,n-1}(\lambda)
+\mu_{n+1}\pi_{j,n+1}(\lambda),\ \ n\geq 2.
\end{cases}
\end{eqnarray}
Let
\begin{eqnarray}\label{eq3-20}
\phi_{j,n}(\lambda)=D_j(\lambda)\pi_{j,n}(\lambda)+F_j(\lambda)\pi_{0,n}(\lambda)
+G_j(\lambda)\pi_{1,n}(\lambda).
\end{eqnarray}
Substitute (\ref{eq3-20}) into the last equality of (\ref{eq3-18}), we have
\begin{eqnarray}\label{eq3-21}
&&D_j(\lambda)[(\lambda+\omega_n+\gamma)\pi_{j,n}(\lambda)-\lambda_{n-1}\pi_{j,n-1}(\lambda)-
\mu_{n+1}\pi_{j,n+1}(\lambda)]-\delta_{j,n}\nonumber\\
&&+F_j(\lambda)[(\lambda+\omega_n+\gamma)\pi_{0,n}(\lambda)-\lambda_{n-1}\pi_{0,n-1}(\lambda)-
\mu_{n+1}\pi_{0,n+1}(\lambda)]\\
&&+G_j(\lambda)[(\lambda+\omega_n+\gamma)\pi_{1,n}(\lambda)-\lambda_{n-1}\pi_{1,n-1}(\lambda)-
\mu_{n+1}\pi_{1,n+1}(\lambda)]=0,\quad n\geq 2.\nonumber
\end{eqnarray}
By the last equalities of (\ref{eq3-12}), (\ref{eq3-13}) and (\ref{eq3-19}), we have
$D_j(\lambda)\delta_{j,n}=\delta_{j,n}$ for $n\geq 2$ and hence $D_j(\lambda)=1$.
\par
Substitute (\ref{eq3-20}) into the first and second equality of (\ref{eq3-18}) and use (\ref{eq3-12}), (\ref{eq3-13}), we have
\begin{eqnarray}\label{eq3-22}
\begin{cases}
(\lambda+\alpha A_{00}(\lambda))F_j(\lambda)
+\alpha A_{10}(\lambda)G_j(\lambda)
=\alpha\lambda\pi_{j,0}(\lambda)-\alpha,\\
\beta A_{01}(\lambda)F_j(\lambda)
+(\lambda+\beta A_{11}(\lambda))G_j(\lambda)
=\beta\lambda\pi_{j,1}(\lambda)-\beta.\\
\end{cases}
\end{eqnarray}
Solving (\ref{eq3-22}) yields (\ref{eq3-8}) and (\ref{eq3-9}).
The proof is complete. \hfill $\Box$
\end{proof}
\par
By Theorem~\ref{th2.1}, we know that \begin{eqnarray*}
\lambda\pi_{j,n}(\lambda)=\lambda\hat{\pi}_{j,n}(\lambda+\gamma)
+\alpha\hat{\pi}_{0,n}(\lambda+\gamma)+\beta\hat{\pi}_{1,n}(\lambda+\gamma).
\end{eqnarray*}
Denote
\begin{eqnarray}\label{eq3-23}
a_n(\lambda)=1-\alpha\hat{\pi}_{0,n}(\lambda+\gamma)
-\beta\hat{\pi}_{1,n}(\lambda+\gamma),\ \ n\geq 0.
\end{eqnarray}
Then, $A_{jn}(\lambda)$ can be represented as
\begin{eqnarray}\label{eq3-24}
A_{jn}(\lambda)=a_{n}(\lambda)-\lambda\hat{\pi}_{j,n}(\lambda+\gamma),
\end{eqnarray}
Hence, by some algebra, $H(\lambda)$ can be represented as
\begin{eqnarray}\label{eq3-25}
H(\lambda)&=&\alpha\beta[a_0(\lambda)\hat{\pi}_{0,1}(\lambda+\gamma)
+a_1(\lambda)\hat{\pi}_{1,0}(\lambda+\gamma)-\lambda \hat{\pi}_{1,0}(\lambda+\gamma)\hat{\pi}_{0,1}(\lambda+\gamma)]\nonumber\\
&&+\alpha a_0(\lambda)\beta(\lambda)+\beta a_1(\lambda)\alpha(\lambda)+\lambda \alpha(\lambda)\beta(\lambda),
\end{eqnarray}
where $\alpha(\lambda)=1-\alpha\hat{\pi}_{0,0}(\lambda+\gamma)$, $\beta(\lambda)=1-\beta\hat{\pi}_{1,1}(\lambda+\gamma)$.
Indeed,
\begin{eqnarray*}
\lambda H(\lambda)&=&(\alpha a_0(\lambda)+\lambda\alpha(\lambda))(\beta a_1(\lambda)+\lambda\beta(\lambda))\\
&&-\alpha\beta(a_0(\lambda)-\lambda\hat{\pi}_{1,0}(\lambda+\gamma))
(a_1(\lambda)-\lambda\hat{\pi}_{0,1}(\lambda+\gamma))\\
&=&\alpha\beta a_0(\lambda)a_1(\lambda)+\alpha\lambda a_0(\lambda)\beta(\lambda)\\
&&+\beta\lambda a_1(\lambda)\alpha(\lambda)+\lambda^2\alpha(\lambda)\beta(\lambda)\\
&&-\alpha\beta a_0(\lambda)a_1(\lambda)+\alpha\beta a_0(\lambda)\lambda\hat{\pi}_{0,1}(\lambda+\gamma)+\alpha\beta a_1(\lambda)\lambda\hat{\pi}_{1,0}(\lambda+\gamma)\\
&&-\alpha\beta \lambda^2\hat{\pi}_{1,0}(\lambda+\gamma)\hat{\pi}_{0,1}(\lambda+\gamma)\\
&=&\lambda \alpha\beta[a_0(\lambda)\hat{\pi}_{0,1}(\lambda+\gamma)
+a_1(\lambda)\hat{\pi}_{1,0}(\lambda+\gamma)-\lambda \hat{\pi}_{1,0}(\lambda+\gamma)\hat{\pi}_{0,1}(\lambda+\gamma)]\\
&&+\lambda[\alpha a_0(\lambda)\beta(\lambda)+\beta a_1(\lambda)\alpha(\lambda)+\lambda \alpha(\lambda)\beta(\lambda)]
\end{eqnarray*}
which implies (\ref{eq3-25}).
\par
\begin{theorem}\label{th3.2}
Let $\Phi(\lambda)=(\phi_{j,n}(\lambda):j,n\in\mathbf{S})$ be the $\tilde{Q}$-resolvent and $\hat{\Pi}(\lambda)=(\hat{\pi}_{j,n}(\lambda):j,n\in\mathbf{Z}_+)$ be the $\hat{Q}$-resolvent. Then,
\begin{eqnarray}\label{eq3-26}
\phi_{j,n}(\lambda)
=\hat{\pi}_{j,n}(\lambda+\gamma)
+\frac{U_j(\lambda)\hat{\pi}_{0,n}(\lambda+\gamma)
+V_j(\lambda)\hat{\pi}_{1,n}
(\lambda+\gamma)}{H(\lambda)},\quad j,n\geq 0,
\end{eqnarray}
where
\begin{eqnarray*}
U_j(\lambda)=\alpha(\lambda+\alpha+\beta)\beta(\lambda)\hat{\pi}_{j,0}(\lambda+\gamma)
+\alpha\beta(\lambda+\alpha+\beta)\hat{\pi}_{1,0}(\lambda+\gamma) \hat{\pi}_{j,1}(\lambda+\gamma),
\end{eqnarray*}
and
\begin{eqnarray*}
V_j(\lambda)=\beta(\lambda+\alpha+\beta)\alpha(\lambda)\hat{\pi}_{j,1}(\lambda+\gamma)
+\alpha\beta(\lambda+\alpha+\beta)\hat{\pi}_{0,1}(\lambda+\gamma) \hat{\pi}_{j,0}(\lambda+\gamma).
\end{eqnarray*}
\end{theorem}
\par
\begin{proof}
By (\ref{eq3-3}), (\ref{eq3-4}) and Theorem~\ref{th2.1}, we know that for any $j,n\geq 0$,
\begin{eqnarray*}
A_{jn}(\lambda)&=&1-\lambda \hat{\pi}_{j,n}(\lambda+\gamma)-\alpha\hat{\pi}_{0,n}(\lambda+\gamma)
-\beta\hat{\pi}_{1,n}(\lambda+\gamma)\\
&=&\lambda[\hat{\pi}_{0,n}(\lambda+\gamma)- \hat{\pi}_{j,n}(\lambda+\gamma)]+A_{0n}(\lambda)\\
&=&\lambda[\hat{\pi}_{1,n}(\lambda+\gamma)- \hat{\pi}_{j,n}(\lambda+\gamma)]+A_{1n}(\lambda).
\end{eqnarray*}
Note that the right hand sides of (\ref{eq3-8}) and (\ref{eq3-9}) are well defined, we can define $F_j(\lambda)$ and $G_j(\lambda)$ for $j=0,1$.
Hence, it follows from Theorem~\ref{th3.1} that for any $j\geq 0$,
\begin{eqnarray*}
\lambda H(\lambda)F_j(\lambda)
&=&\alpha\beta A_{10}(\lambda)A_{01}(\lambda)+\alpha\beta\lambda A_{10}(\lambda)[\hat{\pi}_{0,1}(\lambda+\gamma)- \hat{\pi}_{j,1}(\lambda+\gamma)]\\
&&-\alpha(\lambda+\beta A_{11}(\lambda))A_{00}(\lambda)-\alpha\lambda(\lambda+\beta A_{11}(\lambda))[\hat{\pi}_{0,0}(\lambda+\gamma)- \hat{\pi}_{j,0}(\lambda+\gamma)]
\end{eqnarray*}
and
\begin{eqnarray*}
\lambda H(\lambda)G_j(\lambda)
&=&-\beta\lambda A_{01}(\lambda)+\alpha\beta\lambda A_{01}(\lambda)[\hat{\pi}_{0,0}(\lambda+\gamma)- \hat{\pi}_{j,0}(\lambda+\gamma)]\\
&&-\beta\lambda(\lambda+\alpha A_{00}(\lambda))[\hat{\pi}_{0,1}(\lambda+\gamma)- \hat{\pi}_{j,1}(\lambda+\gamma)].
\end{eqnarray*}
Therefore, by some algebra, one can get
\begin{eqnarray}\label{eq3-27}
&&\lambda H(\lambda)[F_j(\lambda)+\frac{\alpha}{\lambda}
(1+F_j(\lambda)+G_j(\lambda))]\nonumber\\
&=&\alpha\lambda(\lambda+\alpha+\beta)(1 -\beta
\hat{\pi}_{1,1}(\lambda+\gamma))\hat{\pi}_{j,0}(\lambda+\gamma)+
\alpha\beta\lambda(\lambda+\alpha+\beta)\hat{\pi}_{1,0}(\lambda+\gamma) \hat{\pi}_{j,1}(\lambda+\gamma)\nonumber\\
&=:&\lambda U_j(\lambda),\quad j\geq 0.
\end{eqnarray}
Similarly,
\begin{eqnarray}\label{eq3-28}
&&\lambda H(\lambda)[G_j(\lambda)+\frac{\beta}{\lambda}(
1+F_j(\lambda)+G_j(\lambda))]\nonumber\\
&=&\beta\lambda(\lambda+\alpha+\beta)(1-\alpha
\hat{\pi}_{0,0}(\lambda+\gamma))\hat{\pi}_{j,1}(\lambda+\gamma)+
\alpha\beta\lambda(\lambda+\alpha+\beta)\hat{\pi}_{0,1}(\lambda+\gamma) \hat{\pi}_{j,0}(\lambda+\gamma)\nonumber\\
&=:&\lambda V_j(\lambda),\quad j\geq 0.
\end{eqnarray}
By Theorems~\ref{th2.1} and ~\ref{th3.1}, for any $j\geq 2, n\geq 0$,
\begin{eqnarray*}
\phi_{j,n}(\lambda)&=&\pi_{j,n}(\lambda)+F_j(\lambda)\pi_{0,n}(\lambda)
+G_j(\lambda)\pi_{1,n}(\lambda)\\
&=&\hat{\pi}_{j,n}(\lambda+\gamma)+\frac{\alpha\hat{\pi}_{0,n}
(\lambda+\gamma)+\beta\hat{\pi}_{1,n}
(\lambda+\gamma)}{\lambda}\\
&&+F_j(\lambda)[\hat{\pi}_{0,n}(\lambda+\gamma)+\frac{\alpha\hat{\pi}_{0,n}
(\lambda+\gamma)+\beta\hat{\pi}_{1,n}
(\lambda+\gamma)}{\lambda}]\\
&&+G_j(\lambda)[\hat{\pi}_{1,n}(\lambda+\gamma)+\frac{\alpha\hat{\pi}_{0,n}
(\lambda+\gamma)+\beta\hat{\pi}_{1,n}
(\lambda+\gamma)}{\lambda}]\\
&=&\hat{\pi}_{j,n}(\lambda+\gamma)+[F_j(\lambda)+\frac{\alpha}{\lambda} (1+F_j(\lambda)+G_j(\lambda))]\cdot\hat{\pi}_{0,n}(\lambda+\gamma)\\
&&+[G_j(\lambda)+\frac{\beta}{\lambda} (1+F_j(\lambda)+G_j(\lambda))]\cdot\hat{\pi}_{1,n}(\lambda+\gamma),
\end{eqnarray*}
where $F_j(\lambda)$ and $G_j(\lambda)$ are given in (\ref{eq3-8}) and (\ref{eq3-9}). By (\ref{eq3-27}) and (\ref{eq3-28}), we know (\ref{eq3-26}) holds for $j\geq 2, n\geq 0$.
\par
As for $j=0$, by (\ref{eq3-5}) and Theorem~\ref{th2.1},
\begin{eqnarray*}
&&\phi_{0,n}(\lambda)\\
&=&\frac{\lambda(\lambda+\beta A_{11}(\lambda))\pi_{0,n}(\lambda)-\beta\lambda A_{01}(\lambda)\pi_{1,n}(\lambda)}
{\lambda H(\lambda)}\\
&=&\frac{(\lambda+\beta A_{11}(\lambda))[(\lambda+\alpha)\hat{\pi}_{0,n}(\lambda+\gamma)
+\beta\hat{\pi}_{1,n}(\lambda+\gamma)]-\beta A_{01}(\lambda)[(\lambda+\beta)\hat{\pi}_{1,n}(\lambda+\gamma)+\alpha\hat{\pi}_{0,n}(\lambda+\gamma)
]}{\lambda H(\lambda)}\\
&=&\frac{(\lambda+\alpha)(\lambda+\beta A_{11}(\lambda))-\alpha\beta A_{01}(\lambda)}
{\lambda H(\lambda)}\hat{\pi}_{0,n}(\lambda+\gamma)+\frac{\beta[\lambda+\beta A_{11}(\lambda)-(\lambda+\beta) A_{01}(\lambda)]}
{\lambda H(\lambda)}\hat{\pi}_{1,n}(\lambda+\gamma)\\
&=&\hat{\pi}_{0,n}(\lambda+\gamma)+\frac{(\lambda+\alpha)(\lambda+\beta A_{11}(\lambda))-\alpha\beta A_{01}(\lambda)-\lambda H(\lambda)}
{\lambda H(\lambda)}\hat{\pi}_{0,n}(\lambda+\gamma)\\
&&+\frac{\beta[\lambda+\beta A_{11}(\lambda)-(\lambda+\beta) A_{01}(\lambda)]}
{\lambda H(\lambda)}\hat{\pi}_{1,n}(\lambda+\gamma).
\end{eqnarray*}
By the definition of $H(\lambda)$,
\begin{eqnarray*}
&&(\lambda+\alpha)(\lambda+\beta A_{11}(\lambda))-\alpha\beta A_{01}(\lambda)-\lambda H(\lambda)\\
&=&(\lambda+\alpha)(\lambda+\beta A_{11}(\lambda))-\alpha\beta A_{01}(\lambda)-(\lambda+\alpha A_{00}(\lambda))(\lambda+\beta A_{11}(\lambda))+\alpha\beta A_{10}(\lambda)A_{01}(\lambda)\\
&=&\alpha(\lambda+\beta A_{11}(\lambda))(1-A_{00}(\lambda))-\alpha\beta A_{01}(\lambda)(1-A_{10}(\lambda)).
\end{eqnarray*}
\par
On the other hand, by some algebra, one can see that
\begin{eqnarray*}
\lambda U_0(\lambda)&=&\lambda H(\lambda)[F_0(\lambda)+\frac{\alpha}{\lambda}(
1+F_0(\lambda)+G_0(\lambda))]\\
&=&\alpha\beta A_{10}(\lambda)A_{01}(\lambda))
-\alpha(\lambda+\beta A_{11}(\lambda))A_{00}(\lambda))+\alpha(\lambda+\beta A_{11}(\lambda))-\alpha\beta A_{01}(\lambda)\\
&=&\alpha(\lambda+\beta A_{11}(\lambda))(1-A_{00}(\lambda))-\alpha\beta A_{01}(\lambda)(1-A_{10}(\lambda)),
\end{eqnarray*}
\begin{eqnarray*}
\lambda V_0(\lambda)&=&\lambda H(\lambda)[G_0(\lambda)+\frac{\beta}{\lambda}(
1+F_0(\lambda)+G_0(\lambda))]\\
&=&\alpha\beta A_{01}(\lambda)A_{00}(\lambda)
-\beta(\lambda+\alpha A_{00}(\lambda))A_{01}(\lambda)+\beta(\lambda+\beta A_{11}(\lambda))-\beta^2 A_{01}(\lambda)\\
&=&\beta[\lambda+\beta A_{11}(\lambda)-(\lambda+\beta) A_{01}(\lambda)].
\end{eqnarray*}
Therefore, (\ref{eq3-26}) holds for $j=0$. By a similar argument, (\ref{eq3-26}) also holds for $j=1$. The proof is complete.\hfill $\Box$
\end{proof}
\par
We now consider the probability distribution of $C_j$ and the related probabilities $P(C_j\leq t, C_{j,0}<C_{j,1})$ and $P(C_j\leq t, C_{j,1}<C_{j,0})$. It is easy to see that $P(C_j\leq t, C_{j,k}<C_{j,1-k})$ is differentiable in $t$ for $k=0,1$. Let $d_{j,k}(t)=\frac{d}{dt}P(C_j\leq t, C_{j,k}<C_{j,1-k})$ for $k=0,1$. Also,
let $\Delta_{j,k}(\lambda)$ denote the Laplace transform of $d_{j,k}(t)$ for $k=0,1$ and $\Delta_{j}(\lambda)$ denote the Laplace transform of $d_{j}(t)$.
\par
\begin{theorem}\label{th3.3}
For any $j\geq 0$, we have
\begin{eqnarray*}
\Delta_{j,0}(\lambda)=\frac{\alpha(\lambda+\beta)(1-\lambda\phi_{j0}(\lambda))-\alpha\beta(1-\lambda\phi_{j,1}(\lambda))}
{\lambda^2+(\alpha+\beta)\lambda},
\end{eqnarray*}
\begin{eqnarray*}
\Delta_{j,1}(\lambda)=\frac{\beta(\lambda+\alpha)(1-\lambda\phi_{j1}(\lambda))-\alpha\beta(1-\lambda\phi_{j,0}(\lambda))}
{\lambda^2+(\alpha+\beta)\lambda}
\end{eqnarray*}
and
\begin{eqnarray*}
\Delta_j(\lambda)=\frac{\alpha(1-\lambda\phi_{j0}(\lambda))+\beta(1-\lambda\phi_{j1}(\lambda))}
{\lambda+\alpha+\beta},
\end{eqnarray*}
where $\phi_{j,0}(\lambda)$ and $\phi_{j,1}(\lambda)$ are given in Theorem~\ref{th3.2}. In particular,
\begin{eqnarray*}
P(C_{j,0}<C_{j,1})=\frac{\alpha[1+\beta(\phi_{j,1}(0)-\phi_{j0}(0))]
}
{\alpha+\beta},
\end{eqnarray*}
\begin{eqnarray*}
P(C_{j,1}<C_{j,0})=\frac{\beta[1+\alpha(\phi_{j0}(0)-
\phi_{j,1}(0))]}
{\alpha+\beta},
\end{eqnarray*}
where $\phi_{j,0}(\lambda)$ and $\phi_{j,1}(\lambda)$ are given by $(\ref{eq3-26})$.
\end{theorem}
\par
\begin{proof}
By the definitions of $\{M_t:t\geq 0\}$ and $\{N_t:t\geq 0\}$, we know that for any $j\geq 0$,
\begin{eqnarray*}
&&P(C_{j,0}\leq t,C_{j,0}<C_{j,1})=\int_0^{t}d_{j,0}(\tau)d\tau=h_{j,-2}(t),\\
&&
P(C_{j,1}\leq t, C_{j,1}<C_{j,0})=\int_0^{t}d_{j,1}(\tau)d\tau=h_{j,-1}(t)
\end{eqnarray*}
and
\begin{eqnarray*}
P(C_{j}\leq t)=\int_0^{t}d_j(\tau)d\tau=h_{j,-2}(t)+h_{j,-1}(t).
\end{eqnarray*}
Therefore,
$d_{j,0}(t)=h'_{j,-2}(t),\ d_{j,1}(t)=h'_{j,-1}(t)$ and
$d_j(t)=h'_{j,-2}(t)+h'_{j,-1}(t)$.
Hence,
\begin{eqnarray*}
\Delta_{j,0}(\lambda)=\lambda\phi_{j,-2}(\lambda),\quad \Delta_{j,1}(\lambda)=\lambda\phi_{j,-1}(\lambda)
\end{eqnarray*}
and
\begin{eqnarray*}
\Delta_{j}(\lambda)&=&\lambda\phi_{j,-2}(\lambda)+\lambda\phi_{j,-1}(\lambda).
\end{eqnarray*}
\par
By (\ref{eq3-2}) of Lemma~\ref{le3.1}, we know that
\begin{eqnarray*}
(\lambda+\alpha)\lambda \phi_{j,-2}(\lambda)+\alpha\lambda\phi_{j,-1}(\lambda)=\alpha (1-\lambda\phi_{j,0}(\lambda))
\end{eqnarray*}
and
\begin{eqnarray*}
\beta\lambda\phi_{j,-2}(\lambda)+(\lambda+\beta)\lambda\phi_{j,-1}(\lambda)=\beta (1-\lambda\phi_{j,1}(\lambda)).
\end{eqnarray*}
Therefore, by the first two equalities of (\ref{eq3-2}),
\begin{eqnarray*}
\Delta_{j,0}(\lambda)&=&\lambda\phi_{j,-2}(\lambda)=
\frac{\alpha[(\lambda+\beta)(1-\lambda\phi_{j0}(\lambda))-\beta(1-\lambda\phi_{j,1}
(\lambda))]}{\lambda^2+(\alpha+\beta)\lambda},
\end{eqnarray*}
\begin{eqnarray*}
\Delta_{j,1}(\lambda)=\lambda\phi_{j,-1}(\lambda)
=\frac{\beta[(\lambda+\alpha)(1-\lambda\phi_{j1}(\lambda))
-\alpha(1-\lambda\phi_{j,0}(\lambda))]}
{\lambda^2+(\alpha+\beta)\lambda}\\
\end{eqnarray*}
and hence
\begin{eqnarray*}
\Delta_j(\lambda)=\frac{\alpha(1-\lambda\phi_{j0}(\lambda))+\beta(1-\lambda\phi_{j1}(\lambda))}
{\lambda+\alpha+\beta}.
\end{eqnarray*}
Note that $P(C_j<\infty)=\Delta_j(0)=1$, the last two assertions hold. The proof is complete. \hfill $\Box$
\end{proof}
\par
We now consider the mathematical expectation and variance of $C_j$.
\begin{theorem}\label{th3.4}
For any $j\geq 0$,
\begin{eqnarray*}
E[C_j]=\frac{1+\alpha\phi_{j,0}(0)+\beta\phi_{j,1}(0)}{\alpha+\beta}
\end{eqnarray*}
and
\begin{eqnarray*}
E[C^2_j]=\frac{2[1+\alpha\phi_{j,0}(0)+\beta\phi_{j,1}(0)-(\alpha+\beta)
(\alpha\phi'_{j,0}(0)+\beta\phi'_{j,1}(0))]}{(\alpha+\beta)^2},
\end{eqnarray*}
where $\phi_{j,0}(\lambda)$ and $\phi_{j,1}(\lambda)$ are given by $(\ref{eq3-26})$.
\end{theorem}
\par
\begin{proof} By Theorem~\ref{th3.3}, we have
\begin{eqnarray*}
(\lambda+\alpha+\beta)\Delta_j(\lambda)=\alpha(1-\lambda\phi_{j,0}(\lambda))+\beta(1-\lambda\phi_{j,1}(\lambda))
\end{eqnarray*}
Differentiating the above equality yields that
\begin{eqnarray}\label{eq3-29}
(\lambda+\alpha+\beta)\Delta'_j(\lambda)+\Delta_j(\lambda)=-\alpha(\lambda\phi_{j,0}(\lambda))'-\beta(\lambda\phi_{j,1}(\lambda))'
\end{eqnarray}
Let $\lambda=0$ and note that $\Delta_j(0)=1$, we have
\begin{eqnarray*}
E[C_j]=-\Delta'_j(0)=\frac{1+\alpha\phi_{j,0}(0)+\beta\phi_{j,1}(0)}{\alpha+\beta}.
\end{eqnarray*}
Differentiating (\ref{eq3-29}) yields that
\begin{eqnarray*}
&&(\lambda+\alpha+\beta)\Delta''_j(\lambda)+2\Delta'_j(\lambda)\\
&=&-\alpha(\lambda\phi_{j,0}(\lambda))''-\beta(\lambda\phi_{j,1}(\lambda))''\\
&=&-\alpha[\lambda\phi''_{j,0}(\lambda)+2\phi'_{j,0}(\lambda)]
-\beta[\lambda\phi''_{j,1}(\lambda)+2\phi'_{j,1}(\lambda)].
\end{eqnarray*}
Let $\lambda=0$ in the above equality yields that
\begin{eqnarray*}
(\alpha+\beta)\Delta''_j(0)+2\Delta'_j(0)=-2\alpha\phi'_{j,0}(0)
-2\beta\phi'_{j,1}(0).
\end{eqnarray*}
Therefore,
\begin{eqnarray*}
E[C_j^2]&=&\Delta''_j(0)\\
&=&\frac{2(-\Delta'_j(0)-\alpha\phi'_{j,0}(0)
-\beta\phi'_{j,1}(0))}{\alpha+\beta}\\
&=&\frac{2[1+\alpha\phi_{j,0}(0)+\beta\phi_{j,1}(0)-(\alpha+\beta)
(\alpha\phi'_{j,0}(0)+\beta\phi'_{j,1}(0))]}{(\alpha+\beta)^2}.
\end{eqnarray*}
The proof is complete. \hfill $\Box$
\end{proof}
\par
Finally, if $\alpha=0$ or $\beta=0$, we get the following result which is due to Di Crescenzo et al~\cite{DGNR2008}.
\par
\begin{corollary}\label{cor3.1}
\rm{(i)}\ If $\beta=0$, then for any $j\geq 0$,
\begin{eqnarray*}
E[C_j]=\frac{1}{\alpha}+\frac{\hat{\pi}_{j,0}(\alpha)}
{1-\alpha\hat{\pi}_{0,0}(\alpha)}
\end{eqnarray*}
and
\begin{eqnarray*}
E[C^2_j]=\frac{2}{\alpha^2}\left(1+\frac{\alpha\hat{\pi}_{j,0}(\alpha)}
{1-\alpha\hat{\pi}_{0,0}(\alpha)}-\frac{\alpha^2\hat{\pi}'_{j,0}(\alpha)}
{1-\alpha\hat{\pi}_{0,0}(\alpha)}-\frac{\alpha^3\hat{\pi}_{j,0}(\alpha)
\hat{\pi}'_{0,0}(\alpha)}
{(1-\alpha\hat{\pi}_{0,0}(\alpha))^2}\right).
\end{eqnarray*}
\par
\rm{(ii)}\ If $\alpha=0$, then for any $j\geq 0$,
\begin{eqnarray*}
E[C_j]=\frac{1}{\beta}+\frac{\hat{\pi}_{j,1}(\beta)}
{1-\beta\hat{\pi}_{1,1}(\beta)}
\end{eqnarray*}
and
\begin{eqnarray*}
E[C^2_j]=\frac{2}{\beta^2}\left(1+\frac{\beta\hat{\pi}_{j,1}(\beta)}
{1-\beta\hat{\pi}_{1,1}(\beta)}-\frac{\beta^2\hat{\pi}'_{j,1}(\beta)}
{1-\beta\hat{\pi}_{1,1}(\beta)}-\frac{\beta^3\hat{\pi}_{j,1}(\beta)
\hat{\pi}'_{1,1}(\beta)}
{(1-\beta\hat{\pi}_{1,1}(\beta))^2}\right).
\end{eqnarray*}
\end{corollary}
\par
\begin{proof}
If $\beta=0$, by Theorem~\ref{th3.2},
\begin{eqnarray*}
\phi_{j,0}(\lambda)
=\hat{\pi}_{j,0}(\lambda+\alpha)
+\frac{\alpha\hat{\pi}_{j,0}(\lambda+\alpha)
\hat{\pi}_{0,0}(\lambda+\alpha)}{1-\alpha\hat{\pi}_{0,0}(\lambda+\alpha)}
=\frac{\hat{\pi}_{j,0}(\lambda+\alpha)}
{1-\alpha\hat{\pi}_{0,0}(\lambda+\alpha)}.
\end{eqnarray*}
Therefore,
\begin{eqnarray*}
\phi_{j,0}(0)
=\frac{\hat{\pi}_{j,0}(\alpha)}
{1-\alpha\hat{\pi}_{0,0}(\alpha)}
\end{eqnarray*}
and
\begin{eqnarray*}
\phi'_{j,0}(0)
=\frac{\hat{\pi}'_{j,0}(\alpha)}
{1-\alpha\hat{\pi}_{0,0}(\alpha)}+\frac{\alpha\hat{\pi}_{j,0}(\alpha)
\hat{\pi}'_{0,0}(\alpha)}
{(1-\alpha\hat{\pi}_{0,0}(\alpha))^2}.
\end{eqnarray*}
Hence, by Theorem~\ref{th3.3},
\begin{eqnarray*}
E[C_j]=\frac{1+\alpha\phi_{j,0}(0)}{\alpha}=\frac{1}{\alpha}+\frac{\hat{\pi}_{j,0}(\alpha)}
{1-\alpha\hat{\pi}_{0,0}(\alpha)}
\end{eqnarray*}
and
\begin{eqnarray*}
E[C^2_j]&=&\frac{2}{\alpha^2}[1+\alpha\phi_{j,0}(0)-\alpha^2\phi'_{j,0}(0)]\\
&=&\frac{2}{\alpha^2}\left(1+\frac{\alpha\hat{\pi}_{j,0}(\alpha)}
{1-\alpha\hat{\pi}_{0,0}(\alpha)}-\frac{\alpha^2\hat{\pi}'_{j,0}(\alpha)}
{1-\alpha\hat{\pi}_{0,0}(\alpha)}-\frac{\alpha^3\hat{\pi}_{j,0}(\alpha)
\hat{\pi}'_{0,0}(\alpha)}
{(1-\alpha\hat{\pi}_{0,0}(\alpha))^2}\right).\\
\end{eqnarray*}
(i) is proved. The proof of (ii) is similar. \hfill $\Box$
\end{proof}
\par
\section*{Acknowledgement}
\par
This work is supported by the National Natural Science Foundation of China (No. 11771452, No. 11971486).

\par

\end{document}